# On the Maximum–Weight Basis Problem

Brahim Chaourar
Al Imam University (IMSIU)
Riyadh, Saudi Arabia
College of Sciences, P. O. Box 90950, Riyadh 11623, Saudi Arabia
Correspondence address: P. O. Box 287574, Riyadh 11323, Saudi Arabia

**Abstract:** Let M to be a matroid defined on a finite set E. A subset L of E is locked in M if L is 2-connected in M, E\L is 2-connected in the dual M*, and min{r(L), r*(E\L)} ≥ 2. In this paper, we prove that the nontrivial facets of the bases polytope of M are described by the locked subsets. We deduce that finding the maximum–weight basis of M is a polynomial problem for matroids with a polynomial number of locked subsets. This class of matroids is closed under 2-sums and contains uniform matroids.

## 1. Introduction:

Sets and their characteristic vectors will be not distinguished. We refer to Oxley (1992) and Schrijver (1986) for missed definitions and notations about, respectively, matroids and polyhedra.

Let M to be a matroid defined on a finite set E. Independents(M), Bases(M) and the function r are, respectively, the class of independent sets, bases and the rank function of M. M*, Bases(M*) and the function r* are, respectively, the dual matroid, the class of cobases and the dual rank function of M.

For any X⊆E, Bases(X) = {B∈Bases(M) such that |B∩X| = r(X)} and Cobases(X) = {B∈Bases(M*) such that |B∩X| = r*(X)}. The polyhedra Q(M) and P(M) are, respectively, the convex hulls of the independent sets and the bases of M.

Suppose that M and M* are both 2-connected. A subset L of E is locked in M if L is 2-connected (or nonseparable) in M, E\L is 2-connected in M* and min{r(L), r*(E\L)} > 1. It is not difficult to see that if L is locked then both L and E\L are closed respectively in M and M* (That is why we call them locked). We denote by Lockeds(M) the class of locked subsets of M and by ℓ(M) the number of locked subsets of M. Given a nonnegative integer k, a matroid is k–locked if there ℓ(M) is O(|E(M)|$^k$). An infinite class of matroids C is k–locked if all its matroids are k-locked. We also say that C is polynomially locked. It is not difficult to see that the class of lockeds subsets of a matroid M is the union of lockeds subsets of the 2-connected components of M.

Given a weight function c∈R$^E$, the maximum–weight basis problem (MWBP) is the following optimization problem:       Maximize {c(B) such that B∈Bases(M)}.

It is not difficult to see that the corresponding maximum–weight independent problem is (polynomially) equivalent to MWBP.

MWBP is polynomial on |E(M)| and θ, where θ is the complexity of the used matroid oracle. Even if we use the approach introduced by Mayhew (2008) by giving the list of bases (for example) in the input, MWBP is polynomial on |Input|. However, as Robinson and Welsh (1980) note, no matter which of the ways to specify a matroid, the size of the input for a matroid problem on an n-element set is $O(2^n)$. It follows that MWBP is not polynomial in its strict sense, that is on |E(M)|. Our main result is that MWBP is polynomial on |E(M)| for polynomially locked classes of matroids.

The remainder of the paper is organized as follows: in section 2, we give all facets of the bases polytope and deduce our main result, then, in section 3, we describe some polynomially locked classes of matroids. In section 4, we give a polynomial algorithm via a new matroid oracle for testing if a given matroid is uniform or not, and finally we conclude in section 5.

**2. The main result:**

A description of Q(M) was given by Edmonds (1971) as follows.

**Theorem**: Q(M) is the set of all $x \in R^E$ such that

$$x(e) \geq 0 \quad \text{for any } e \in E \tag{1}$$
$$x(A) \leq r(A) \quad \text{for any } A \subseteq E \tag{2}$$

Later, a minimal description of Q(M) was given also by Edmonds (Giles 1975) as follows.

**Theorem**: The inequality (2) is a facet of Q(M) if and only if A is closed and 2-connected.

It is not difficult to see that P(M) is the set of all $x \in R^E$ satisfying the inequalities (1), (2) and

$$x(E) = r(E) \tag{3}$$

It seems natural to think that the inequality (2) is a facet of P(M) if and only if A is closed and 2-connected. This is not true because:

**Lemma 1:** If the inequality (2) is a facet of P(M) then A is a locked subset of M.

Proof: It suffices to prove that if X is closed and 2-connected but E\L is not 2-connected in the dual then the inequality (2) is not a facet. In fact, there exist A and B two disjoint subsets of E such that E\X = A∪B and r*(E\X) = r*(A)+r*(B), that is, |E\X|–r(E)+r(X) = |A|–r(E)+r(E\A)+|B|–r(E)+r(E\B). It follows that: r(E)+r(X) = r(E\A)+r(E\B) ≥ x(E\A)+x(E\B) = x(E)+x(X), which implies the inequality (2). So the inequality (2) is redundant and cannot be a facet. □

We give now a minimal description of P(M). A part of the proof is inspired from a proof given by Pulleyblank (1989) to describe the nontrivial facets of Q(M). Independently, Fujishige (1984), and Feichtner and Sturmfels (2005), gave a characterization of nontrivial facets of P(M). We give here below a new and complete formulation with a new proof.

**Theorem 2:** A minimal description of P(M) is the set of all $x \in R^E$ satisfying the equality (3) and the following inequalities:

$$x(P) \leq 1 \quad \text{for any parallel closure P of E} \tag{4}$$
$$x(S) \geq |S|-1 \quad \text{for any coparallel closure S of E} \tag{5}$$
$$x(L) \leq r(L) \quad \text{for any locked subset L of E} \tag{6}$$

Proof:

Without loss of generality, we can suppose that M is without parallel or coparallel elements so the inequalities (5) become as (1) and the inequalities (4) become as follows:

$$x(e) \leq 1 \quad \text{for any } e \in E \tag{7}$$

Let C(M) be the cone generated by the incidence vectors of the bases of M. It suffices to prove that the minimal description of C(M) is given by (1) and the following inequalities:

$$x(e) \leq x(E)/r(E) \quad \text{for any } e \in E \tag{8}$$
$$x(L)/r(L) \leq x(E)/r(E) \quad \text{for any locked subset L of E} \tag{9}$$

It is not difficult to see via induction and operations of deletion and contraction that the inequalities (1) and (8) are facets of C(M). It remains to prove that the inequality (9) is a facet of C(M) if and only if L is a locked subset of M. According to Lemma 1, it suffices to prove the inverse way.

Note that (9) is equivalent to the following inequality:

$$(r(L)-r(E)) \, x(L) + r(L) \, x(E \backslash L) \geq 0 \quad \text{for any locked subset L of E} \tag{10}$$

Let $a\,x \geq 0$ be a valid inequality for C(M) which is tight for all $B \in Bases(L)$.

**Claim 1:** $a_j = a_k$ for all j and k of L.

Suppose that is not true. Let X = {j∈L such that $a_j$ takes minimum value over L}, Y = L\X and B∈Bases(L)∩Bases(Y).

Since L is 2-connected in M, and since, by assumption, X is a strict subset of L, then r(X) > |B∩X|. Thus it exists e∈X\B such that (B∩X)∪{e} is an independent set of M. It follows that there exists f∈B∩Y such that B' = B\{f})∪{e}∈Bases(L). But:
a(B') = a(B)–a(f)+a(e) < a(B) a contradiction.

**Claim 2:** For any X⊆E, B∈Bases(X) if and only if E\B∈Cobases(E\X).

It suffices to prove one way and use duality for the other way.

Let B∈Bases(X) then
|B∩X| = r(X) = |X|–r*(E)+r*(E\X) = |E|–|E\X|–r*(E)+r*(E\X) = r(E)–|E\X|+r*(E\X).

Thus,
|(E\B)∩(E\X)| = |E\X|–|B∩(E\X)| = |E\X|–|B|+|B∩X| = |E\X|–|B|+r(E)–|E\X|+r*(E\X) = r*(E\X).
Since E\B is a basis in the dual, then E\B∈Cobases(E\X).

**Claim 3:** $a_j = a_k$ for all j and k of E\L.

Using claim 2, E\L is 2-connected and a similar argument on E\B as in claim 1, we conclude.

**Claim 4:** $a\,x \geq 0$ is a multiple of inequality (10).

By claims 1 and 3, $a\,x \geq 0$ becomes: $a_L \, x(L) + a_{E \backslash L} \, x(E \backslash L) \geq 0$. Thus, for B∈Bases(L), we have:
$a_L \, |B \cap L| + a_{E \backslash L} \, |B \cap (E \backslash L)| = 0$, that is, $a_L \, r(L) + a_{E \backslash L} \, (r(E)-r(L)) = 0$. But $a_L = r(L)-r(E)$ and $a_{E \backslash L} = r(L)$ is a solution of this equation and we conclude. □

The main result follows.

**Corollary 3:** Let C to be a polynomially locked class of matroids defined on a finite set E and $c \in R^E$. Then there exists a polynomial algorithm on the size of E for solving MWBP in any matroid M of C.

Proof:

Let M be a matroid of C. M can be described by its locked subsets in the input of MWBP. MWBP is equivalent for optimizing on P(M), which is also equivalent to separating on P(M). Since the number of facets of P(M) is bounded by 2|E|+ℓ(M) then separating can be done on O(|E|+ℓ(M)). But M is a member of a polynomially locked class of matroids, then ℓ(M) is bounded by a polynomial on the size of E. □

### 3. Some polynomially locked classes of matroids:

We need some results on matroids.

The following results were proved or can be deduced from results proved by Chaourar (2008 and 2011).

**Theorem 4:**
  (i)   Uniform matroids are the matroids with no locked subsets.
  (ii)  Series parallel extensions preserve k-lockdness for $k \geq 1$.
  (iii) 2-sums preserve k-lockdness for $k \geq 1$.

The following result was proved by Chaourar and Oxley (2003).

**Theorem 5:** The excluded minors for the class of matroids that can be constructed from uniform matroids by a sequence of series extensions, parallel extensions, or direct sums are $M(K_4), W^3, P_6, Q_6$, and $U_{2,4} \oplus_2 U_{2,4}$.

And finally the following result should be credited for Walton (1981) (or see Chaourar 2011).

**Theorem 6:** The excluded minors for the class of matroids that can be constructed from uniform matroids by a sequence of 1- sums and 2-sums are $M(K_4), W^3, P_6$ and $Q_6$.

It follows the following direct results combining Corollary 3, Theorems 4, 5 and 6.

**Corollary 7:** MWBP is polynomial in the size of the ground set (the class of locked subsets of the matroid describes it in the input) for the following classes of matroids:

  (1) Matroids without $M(K_4), W^3, P_6, Q_6$, and $U_{2,4} \oplus_2 U_{2,4}$ as minors;
  (2) Matroids without $M(K_4), W^3, P_6$ and $Q_6$ as minors.

### 4. Recognizing uniform matroids:

Given a nonnegative integer k, a natural k–locked oracle is as follows:

Instance: A nonnegative integer k and a matroid M.

Question: Is M k–locked? And if the answer is Yes then give the list of locked subsets of M.

Note that the running time complexity of this oracle is $O(|E|^k)$ because we have to count at most $|E|^k+1$ elements of Lockeds(M) to know the answer.

It follows that testing if a given matroid is uniform can be done in polynomial running time on the size of the ground set by using the 0–locked oracle and Theorem 4.

### 5. Conclusion:

We have proved that MWBP is polynomial on the size of the ground set for polynomially locked matroids and, in particular, for matroids without $M(K_4), W^3, P_6$ and $Q_6$ as minors.

We have also given a polynomial time algorithm on the size of the ground set and the complexity of the 0–locked oracle for testing if a matroid is uniform or not.

Further investigations can be done to find more larger classes of polynomially locked matroids.